\documentclass[12pt,a4paper,leqno]{article}
\usepackage{etoolbox} \usepackage{comment}
\usepackage{amsfonts}
\usepackage[margin=1in]{geometry}
\usepackage{times}
\usepackage{amsthm,amssymb}
\usepackage{amsmath,xypic}
\usepackage{pdflscape}
\usepackage[export]{adjustbox} 
\usepackage[cal=boondoxo,calscaled=.96,scr=rsfs]{mathalpha} \usepackage{multicol}
\usepackage{hyphenat}
\usepackage[usenames,dvipsnames]{color}
\usepackage{rotating}
\usepackage{longtable}
\usepackage{caption}
\usepackage{tcolorbox}

\usepackage{epigraph}
\usepackage{enumitem}
\setlist[enumerate]{listparindent=0.5in}
\newcommand{\be}{\begin{equation}}
\newcommand{\ee}{\end{equation}}
\newcommand{\bes}{\begin{equation*}}
\newcommand{\ees}{\end{equation*}}
\newcommand{\bea}{\begin{eqnarray}}
\newcommand{\eea}{\end{eqnarray}}
\newcommand{\beas}{\begin{eqnarray}}
\newcommand{\eeas}{\end{eqnarray}}
\newcommand{\ben}{\begin{note}}
\newcommand{\een}{\end{note}}
\newcommand{\bexl}{\vskip0.1em\noindent\hrulefill\vskip1em\begin{ExerciseList}}
\newcommand{\eexl}{\end{ExerciseList}\hrulefill}

\newcommand{\bthm}{\begin{theorem}}
\newcommand{\ethm}{\end{theorem}}
\newcommand{\bpro}{\begin{prop}}
\newcommand{\epro}{\end{prop}}
\newcommand{\bcor}{\begin{corollary}}
\newcommand{\ecor}{\end{corollary}}
\newcommand{\bcon}{\begin{conjecture}}
\newcommand{\econ}{\end{conjecture}}
\newcommand{\bp}{\begin{proof}}
\newcommand{\ep}{\end{proof}}
\newcommand{\blem}{\begin{lemma}}
\newcommand{\elem}{\end{lemma}}
\newcommand{\bn}{\begin{note}}
\newcommand{\en}{\end{note}}
\newcommand{\benum}{\begin{enumerate}}
\newcommand{\eenum}{\end{enumerate}}
\newcommand{\bed}{\begin{defn}}
\newcommand{\eed}{\end{defn}}
\newcommand{\brem}{\begin{remark}}
\newcommand{\erem}{\end{remark}}

\newcommand{\btik}{\begin{tikzpicture}\begin{axis}[scale=0.5,axis y line=center, axis x line=middle]}
\newcommand{\etik}{\end{axis}\end{tikzpicture}}

\let\into=\hookrightarrow
\let\mapsto=\longmapsto

\newcommand{\upperRomannumeral}[1]{\uppercase\expandafter{\romannumeral#1}}

 \usepackage{tikz}
\usetikzlibrary{mindmap,backgrounds}
\usetikzlibrary{automata}
\usepackage{graphicx}
\usepackage[pagewise]{lineno}
\usepackage{stackengine}
\usepackage{stmaryrd}

\newcommand{\authornamemode}{\usepackage[backend=biber,style=authoryear]{biblatex}\setlength\bibitemsep{\baselineskip}}
\usepackage{titlesec}
\newtoggle{arxiv}
\toggletrue{arxiv}
\iftoggle{arxiv}{
	\usepackage{url}
\usepackage{natbib}
	\bibliographystyle{plainnat}
	\let\cite=\citep
	\typeout{In the arxiv mode}
	\usepackage{fancyhdr}
	\usepackage[colorlinks,citecolor=blue]{hyperref}
}
{
	\usepackage{git}
    \usepackage[pagewise]{lineno}
\authornamemode
	\usepackage{fancyhdr}
	\usepackage[colorlinks,citecolor=blue,colorlinks=true,hyperindex, citecolor=blue, urlcolor=blue]{hyperref}
	\let\cite=\parencite
\addbibresource{../../master/master6.bib}
\addbibresource{hoshi-bib.bib}
\addbibresource{mochizuki-bib.bib}
\addbibresource{uchida-bib.bib}
\addbibresource{mochizuki-flowchart.bib}
}
\newtoggle{draft}
\togglefalse{draft}

\newcommand{\preliminary}{{\\ \normalsize \textcolor{blue}{Preliminary version for comments}}{\relax}}
\newtheorem{theorem}[equation]{Theorem}      \newtheorem{lemma}[equation]{Lemma}          \newtheorem{corollary}[equation]{Corollary}  \newtheorem{proposition}[equation]{Proposition}

\theoremstyle{definition}
\newtheorem{conj}[equation]{Conjecture}

\theoremstyle{definition}
\newtheorem{defn}[equation]{Definition}
\theoremstyle{remark}

\theoremstyle{definition}
\newtheorem{remark}[equation]{Remark}

\numberwithin{equation}{section}

\newcommand{\para}{\subsection{}}
\titleformat{\subsection}[runin]{\normalfont\bfseries}{\S\ \thesubsection}{.5em}{}[{\ \ }]
\titlespacing{\subsection}{0pt}{1.5ex plus .1ex minus .2ex}{0pt}

\let\into=\hookrightarrow
\let\isom=\simeq

\let\tensor=\otimes
\newcommand{\A}{\mathscr{A}}
\newcommand{\abs}[1]{\left\vert#1\right\vert}

\newcommand{\bF}{{\bar{F}}}
\newcommand{\bQ}{{\bar{\Q}}}

\newcommand{\C}{{\mathbb C}}

\newcommand{\F}{{\mathbb F}}

\newcommand{\gal}{{\rm Gal}}

\newcommand{\N}{\mathscr{N}}

\newcommand{\Q}{{\mathbb Q}}
\newcommand{\R}{{\mathbb R}}

\newcommand{\Z}{{\mathbb Z}}

\renewcommand{\O}{{\mathscr O}}
\renewcommand{\P}{{\mathbb P}}
\renewcommand{\wp}{{\mathfrak p}}

\newcommand{\fm}{{\mathfrak{M}}}

\let\fm=\fa
\newcommand{\mapright}[1]{{\xymatrix{{}\ar[r]^{#1}&{}}}}

\renewcommand{\bpro}{\begin{proposition}}
	\renewcommand{\epro}{\end{proposition}}
\renewcommand{\bcon}{\begin{conj}}
	\renewcommand{\econ}{\end{conj}}

\title{Construction of  Arithmetic Teichmuller spaces II: Towards Diophantine Estimates
\preliminary
}
\author{Kirti Joshi}

\newcommand{\Address}{\bigskip\noindent{\footnotesize\textsc{{Math. department, University of Arizona, 617 N Santa Rita, Tucson
		85721-0089, USA.}}\par\nopagebreak 
\noindent\textit{Email:}	\texttt{kirti@math.arizona.edu}}}

\begin{document}
	\maketitle
\epigraphwidth0.55\textwidth

\begin{abstract}
	This paper deals with three consequences of the existence of Arithmetic Teichmuller spaces of \url{https://arxiv.org/abs/2106.11452}.  Let  $\mathscr{X}_{F,\mathbb{Q}_p}$ (resp. $B=B_{\mathbb{Q}_p}$) be the complete Fargues-Fontaine curve (resp. the ring)  constructed by Fargues-Fontaine with the datum $F={\mathbb{C}_p^\flat}$ (the tilt of $\mathbb{C}_p$), $E=\mathbb{Q}_p$. Fix an odd prime $\ell$, let $\ell^*=\frac{\ell-1}{2}$. 
	\begin{enumerate}
		\item The construction (\S 7) of an uncountable subset  $\Sigma_{F}\subset \mathscr{X}_{F,\mathbb{Q}_p}^{\ell^*}$ with  a simultaneous valuation scaling property (Theorem 7.8.1), Galois action and other symmetries. 
		\item Now fix a Tate elliptic curve over a finite extension of $\mathbb{Q}_p$.  The  existence of $\Sigma_{F}$ leads  to the construction (\S 9) of a  set $\widetilde{\Theta}\subset B^{\ell^*}$  consisting of lifts (to $B$),  of values  (lying in different untilts provided by $\Sigma_{F}$) of a chosen theta-function  evaluated at $2\ell$-torsion points on  the chosen  elliptic curve.  The construction of $\widetilde{\Theta}$ can be easily adelized.  Moreover I also prove a lower bound (Theorem 10.1.1) for the size of $\widetilde{\Theta}$ (here size is defined in terms of the Fr\'echet structure of $B$). 
		\item   I also demonstrate  (in \S 11) the existence of ``log-links'' in the theory of [Joshi 2021]. 
	\end{enumerate}
\end{abstract}

\lhead{}

\iftoggle{draft}{\pagewiselinenumbers}{\relax}
\newcommand{\act}{\curvearrowright}
\newcommand{\lmp}{{\Pi\act\Ot}}
\newcommand{\lmpi}{{\lmp}_{\int}}
\newcommand{\lmpf}{\lmp_F}
\newcommand{\Om}{\O^{\times\mu}}
\newcommand{\Omf}{\O^{\times\mu}_{\bF}}
\renewcommand{\N}{\mathbb{N}}
\newcommand{\yoga}{Yoga}
\newcommand{\gl}[1]{{\rm GL}(#1)}
\newcommand{\bK}{\overline{K}}
\newcommand{\reptrip}{\rho:G_K\to\gl V}
\newcommand{\reptripp}[1]{\rho\circ\alpha:G_{\ifstrempty{#1}{K}{{#1}}}\to\gl V}
\newcommand{\benumlab}{\begin{enumerate}[label={{\bf(\arabic{*})}}]}
\newcommand{\ord}{\mathop{\rm ord}\nolimits}	
\newcommand{\kcs}{K^\circledast}
\newcommand{\lcs}{L^\circledast}
\renewcommand{\A}{\mathbb{A}}
\newcommand{\bfq}{\bar{\mathbb{F}}_q}
\newcommand{\tripod}{\P^1-\{0,1728,\infty\}}

\newcommand{\vseq}[2]{{#1}_1,\ldots,{#1}_{#2}}
\newcommand{\anab}[4]{\left({#1},\{#3 \}\right)\anabelmap\left({#2},\{#4 \}\right)}

\newcommand{\gln}{{\rm GL}_n}
\newcommand{\glo}[1]{{\rm GL}_1(#1)}
\newcommand{\glt}[1]{{\rm GL_2}(#1)}

\newcommand{\iut}{\cite{mochizuki-iut1, mochizuki-iut2, mochizuki-iut3,mochizuki-iut4}}
\newcommand{\topics}{\cite{mochizuki-topics1,mochizuki-topics2,mochizuki-topics3}}

\newcommand{\linv}{\mathfrak{L}}
\newcommand{\bedef}{\begin{defn}}
\newcommand{\eedef}{\end{defn}}
\renewcommand{\act}[1][]{\overset{#1}{\curvearrowright}}
\newcommand{\bfx}{\overline{F(X)}}
\newcommand{\anabelmap}{\leftrightsquigarrow}
\newcommand{\ban}[1][G]{\mathscr{B}({#1})}
\newcommand{\pit}{\Pi^{temp}}
 
 \newcommand{\bL}{\overline{L}}
 \newcommand{\bkm}{\bK_M}
 \newcommand{\vbk}{v_{\bK}}
 \newcommand{\vbkm}{v_{\bkm}}
\newcommand{\ocs}{\O^\circledast}
\newcommand{\ot}{\O^\triangleright}
\newcommand{\ocsk}{\ocs_K}
\newcommand{\otk}{\ot_K}
\newcommand{\ok}{\O_K}
\newcommand{\oko}{\O_K^1}
\newcommand{\oks}{\ok^*}
\newcommand{\Qpb}{\overline{\Q}_p}
\newcommand{\Qpbh}{\widehat{\overline{\Q}}_p}
\newcommand{\tr}{\triangleright}
\newcommand{\ocpt}{\O_{\C_p}^\tr}
\newcommand{\ocpf}{\O_{\C_p}^\flat}
\newcommand{\sG}{\mathscr{G}}
\newcommand{\sxfe}{\mathscr{X}_{F,E}}
\newcommand{\sxfep}{\mathscr{X}_{F,E'}}
\newcommand{\loglt}{\log_{\sG}}
\newcommand{\fc}{\mathfrak{t}}
\newcommand{\ku}{K_u}
\newcommand{\kup}{\ku'}
\newcommand{\kt}{\tilde{K}}
\newcommand{\sGpf}{\sG(\O_K)^{pf}}
\newcommand{\hgm}{\widehat{\mathbb{G}}_m}
\newcommand{\bE}{\overline{E}}
\newcommand{\sY}{\mathscr{Y}}
\newcommand{\syfe}{\mathscr{Y}_{F,E}}
\newcommand{\syfqp}[1]{\mathscr{Y}_{\cptl{#1},\Q_p}}
\newcommand{\syfqpe}[1]{\mathscr{Y}_{{#1},E}}
\newcommand{\fJ}{\mathfrak{J}}
\newcommand{\fM}{\mathfrak{M}}
\newcommand{\locvar}{local arithmetic-geometric anabelian variation of fundamental group of $X/E$ at $\wp$}
\newcommand{\fjxep}{\fJ(X,E,\wp)}
\newcommand{\fjxe}{\fJ(X,E)}
\newcommand{\fpc}[1]{\widehat{{\overline{\F_p(({#1}))}}}}
\newcommand{\cpt}{\C_p^\flat}
\newcommand{\cptl}[1]{\C_{p,{#1}}^\flat}
\newcommand{\fja}[1]{\fJ^{\rm arith}({#1})}
\newcommand{\ainfe}{A_{\inf,E}(\O_F)}
\renewcommand{\ainfe}{W_{\O_E}(\O_F)}
\newcommand{\gmh}{\widehat{\mathbb{G}}_m}
\newcommand{\sE}{\mathscr{E}}
\newcommand{\bpi}{B^{\varphi=\pi}}
\newcommand{\bpip}{B^{\varphi=p}}
\newcommand{\onto}{\twoheadrightarrow}

\newcommand{\cpmax}{{\C_p^{\rm max}}}
\newcommand{\xan}{X^{an}}
\newcommand{\yan}{Y^{an}}
\newcommand{\bPi}{\overline{\Pi}}
\newcommand{\bPit}{\bPi^{\rm{\scriptscriptstyle temp}}}
\newcommand{\Pit}{\Pi^{\rm{\scriptscriptstyle temp}}}
\renewcommand{\pit}[1]{\Pi^{\scriptscriptstyle temp}_{#1}}
\newcommand{\pitk}[2]{\Pi^{\scriptscriptstyle temp}_{#1;#2}}
\newcommand{\pio}[1]{\pi_1({#1})}
\newcommand{\fTeich}{\widetilde{\fJ(X/L)}}
\newcommand{\ssep}{\S\,} \newcommand{\vphi}{\varphi}
\newcommand{\sgt}{\widetilde{\sG}}

\setcounter{tocdepth}{1}

\tableofcontents

\togglefalse{draft}
\newcommand{\FF}{\cite{fargues-fontaine}}
\iftoggle{draft}{\pagewiselinenumbers}{\relax}

\newcommand{\attportion}{Sections~\ref{se:number-field-case}, \ref{se:construct-att}, \ref{se:relation-to-iut}, \ref{se:self-similarity} and \ref{se:applications-elliptic}}

\newcommand{\four}{Sections~\ref{se:grothendieck-conj}, \ref{se:untilts-of-Pi}, and \ref{se:riemann-surfaces}}

\numberwithin{equation}{subsection}
\newcommand{\tcp}{\widetilde{\C}_p}
\newcommand{\tK}{\widetilde{K}}
\newcommand{\tF}{\widetilde{F}}
\newcommand{\bfp}{\overline{\F}_p}
\newcommand{\sxqp}{\mathscr{X}_{\cpt,\Q_p}}
\newcommand{\syQp}{\mathscr{Y}_{\cpt,\Q_p}}
\newcommand{\sxQp}{\mathscr{X}_{\cpt,\Q_p}}
\newcommand{\ttxlv}{\tilde{\Theta}_{X,\ell;v}}
\newcommand{\tm}[1]{\theta_M({#1})}
\renewcommand{\fm}{\frak{m}}
\newcommand{\ells}{{\ell^*}}
\newcommand{\cpmaxt}{\C_p^{max\flat}}
\newcommand{\sL}{\mathscr{L}}

\numberwithin{equation}{subsection}
\newcommand{\ttheta}{\widetilde{\Theta}}
\section{Introduction} \para The main result of this paper is Theorem~\ref{th:main}. Let me provide some motivation for the constructions of this paper which lead to the main theorem. Let  $q$ be  the Tate parameter of a Tate elliptic curve over a $p$-adic field $E$ (a $p$-adic field is a finite extension of $\Q_p$ for some prime $p$). Mochizuki's \iut\ rests on the remarkable assertion that at some level in his theory (\iut),  one can set (intentionally written in quotes--for this discussion) $$\text{``}q=q^{j^2}\text{''}$$ for integers $j$ and even set
\be\label{eq:ansatz} \text{``}q=(q,q^{2^2},\cdots, q^{\ells^2})\text{''}
\ee (here and elsewhere in this paper $\ells=\frac{\ell-1}{2}$ for a prime $\ell\geq 3$ which will be chosen later). This assertion is (roughly speaking) the theory of $\theta$-links of \iut\ and it, together with Mochizuki's Yoga of Indeterminacies (described in loc. cit.) leads to the construction of a fundamental set of  interest  which plays a central role in the Diophantine estimates of \iut.  I shall refer to \eqref{eq:ansatz} as \emph{Mochizuki's ansatz}. The important property of Mochizuki's  ansatz is that it underlies \emph{a simultaneous scaling of valuations} (about which I will say more in \ssep\ref{pa:ansatz-explained-in-intro}).

\para Using the notion of  Arithmetic Teichmuller space $\fjxe$ and especially $\fjxe_{\cpt}$ developed in  \cite[Section 8]{joshi-tempered}, here  (see \ssep\ref{se:ansatz}) I provide a natural and intrinsic construction of Mochizuki's ansatz.   In fact, in (\ssep \ref{se:ansatz}), I construct  a set $\Sigma_F$, which I call \emph{Mochizuki's Ansatz}, each point of which provides an equation \eqref{eq:ansatz} and show that Mochizuki's Ansatz $\Sigma_F$ leads to an intrinsic and independent construction  of  a set $\ttheta$ (\ssep\ref{se:construction-ttheta}) similar to the set constructed in \cite[Corollary 3.12]{mochizuki-iut3}. Here $\ttheta$ is constructed  explicitly as a subset   $$\ttheta\subset B^\ells=B\times B\times \cdots\times B$$
where $B=B_{\Q_p}$ is the ring constructed in \cite[Chap. 2]{fargues-fontaine}.  

\para\label{pa:ansatz-explained-in-intro} From my point of view, Mochizuki's ansatz \eqref{eq:ansatz}  has an elegant  formulation (\ssep\ref{se:ansatz}):
\be\label{eq:ansatz2} p=([p^{\flat}],[(p^{\flat})^{2^2}],
\cdots,[(p^{\flat})^{\ells^2}]) \text{ that is } p=[p^\flat]=[p^\flat]^{2^2}=\cdots =[p^\flat]^{\ells^2},
\ee
where $([p^\flat]-p)\subset W(\O_{\cpt})$ is the prime ideal giving rise to some characteristic zero untilt $K$ of $\cpt$ and where $\cpt$ is the tilt of $\C_p$--note that $\cpt$ is an algebraically closed perfectoid (hence complete) field of characteristic $p>0$ (see \cite[Lemma 3.4]{scholze}). This makes the construction and the existence  of $\Sigma_F$ abundantly clear--it is the subset of closed  points of degree one $(y_1,\ldots,y_\ells)\in\sxQp^\ells$  where $[p^\flat]-p$ corresponds to $y_1$ and the relationship \eqref{eq:ansatz2} holds among the tilts of $p$ corresponding to each $y_i$ i.e. $([{p^\flat}^{i^2}]-p)$ is the maximal ideal of $W(\O_{\cpt})$ corresponding to $y_i$. In particular, Mochizuki's ansatz \eqref{eq:ansatz2} exists independently of  Mochizuki's Anabelian Reconstruction Theory (\topics, \iut) and in greater generality--and hence I have chosen to call it by a different name. 

That Mochizuki's ansatz \eqref{eq:ansatz2} also provides the crucial, and non-trivial, simultaneous valuation scaling property between valuations of the respective residue fields of $(y_1,\ldots,y_\ells)$ is established in Theorem~\ref{pr:lift-vals}. 

Now fix a $p$-adic theta-function on a Tate elliptic curve over a finite extension of $\Q_p$. In \ssep\ref{pa:lifting-values-to-B}, I  explain how one may lift theta-values at $2\ell$-torsion points to $B$ and how non-uniqueness of lifts provides additional degree of freedom in making estimates (the idea of lifting to $B$ is detailed in \cite{joshi-tempered}). As was noted in \cite{joshi-tempered}, the idea of lifting to $B$ is analogous to Mochizuki's idea of working with $\theta$-value-monoids. This together with Mochizuki's Ansatz $\Sigma_F$ leads one naturally to the construction of the subset $\ttheta$ mentioned above. In Theorem~\ref{thm:theta-pilot-object-appears}, I show that Mochizuki's ``\emph{Theta-pilot object}'' appears quite naturally from the point of view this paper and  the Arithmetic Teichmuller Theory of \cite{joshi-tempered} and it is indeed subject to a meaningful variation--this is the existence of the set $\ttheta$. In Theorem~\ref{th:main}, I prove a lower bound on $\ttheta$ of the type established in \cite[Corollary 3.12]{mochizuki-iut3}.  

\para \emph{In a very broad sense, the method of construction of $\ttheta$ carried out here is inspired by Mochizuki's ideas and follows Mochizuki's rubric for the construction of his set, but there are many key innovations (such as $\Sigma_F$) introduced here and in \cite{joshi-tempered}}.

\para The point of Mochizuki's ansatz  is that it asserts that Arithmetic Teichmuller Theory (of \cite{joshi-tempered})  admits
\benumlab
\item sets such as $\Sigma_F$ providing a simultaneous scaling of valuations (Theorem \ref{pr:lift-vals});
\item and  points of $\Sigma_F$  typically correspond to distinct Berkovich geometries (this aspect is detailed in \cite{joshi-tempered})--so this is not a mere rescaling of valuations; and finally that
\item   one can  ``glue''  the  data of arithmetic line bundles arising from such scaling of valuations i.e. changes in Berkovich geometries (this is the existence of the set $\ttheta$).
\eenum
\brem 
In \iut, Mochizuki does not consider changes in Berkovich geometries but asserts gluing of different scheme theories.
\erem

\para   For an introduction to the set constructed by Mochizuki's \cite[Corollary 3.12]{mochizuki-iut3}  and how a bound on its size  has deep Diophantine consequences, readers may also see other discussions of \cite[Corollary 3.12]{mochizuki-iut3}  for example \cite{fucheng}, \cite{yamashita},  \cite{dupuy2020statement}. Note that in \cite{dupuy2020statement}, Mochizuki's Corollary 3.12 is stated as Conjecture 1.0.1.

\para Mochizuki's approach suggests \emph{very optimistically} that an assertion \emph{similar}  to \cite[Corollary 3.12]{mochizuki-iut3} \emph{might} be expected to hold for the subset $\ttheta$ constructed here (see \ssep\ref{pa:bound-of-cor312}) and indeed, such a  bound is established in Theorem~\ref{th:main}. But \emph{to be absolutely and  perfectly clear} and  as is clarified below,  Theorem~\ref{th:main}  is structurally similar \emph{but not the same} as \cite[Corollary 3.12]{mochizuki-iut3}.  Importantly,  Theorem~\ref{th:main} is still local (i.e. takes place at a single prime),  but it can be globalized, i.e. adelized, quite easily in the presence of a number field using the adelic arithmetic Teichmuller spaces of \cite{joshi-teich}--but this adelic version is not treated here.

\para  In \cite{mochizuki-iut4}, Mochizuki shows that for his (global  i.e. adelic) version of  the set $\ttheta$ there is an upper bound of the form (my notation here is multiplicative while Mochizuki's notation is additive using the standard (real-valued) logarithm) 
\be\label{eq:mochizukis-inequality} 
\prod_v \abs{q_v^{1/2\ell}}^\ells_{\C_p}\leq { |\ttheta|_{adelic}} \leq A(\ell)\left(\prod_v \abs{q_v^{1/2\ell}}_{\C_p}\right)^{B(\ell)}
\ee where the product is now over all primes (at primes of good reduction one has $q_v=1$), and calculates  $A(\ell), B(\ell)$ explicitly showing that $A(\ell)$ depends on $\ell$, the conductor of the underlying elliptic curve and its field of definition while $B(\ell)\geq \ells$ is independent of the curve but depends on $\ell$)  and that this explicit inequality \eqref{eq:mochizukis-inequality} leads to Diophantine applications such as the Szpiro's conjecture (see \iut, \cite{mochizuki-gaussian}) and also for gentler introductions to these ideas  see \cite{fesenko-iut},  \cite{fucheng}, \cite{yamashita}, \cite{dupuy2020probabilistic}).

\para \emph{To be perfectly clear:}
\benumlab
\item No Diophantine consequences of Theorem~\ref{th:main} are being claimed in this paper. 
\item Notably I do not claim to prove that $\ttheta$ constructed here satisfies the upper bound of \eqref{eq:mochizukis-inequality}.
\item Results of this paper do not give a new proof of \cite[Corollary 3.12]{mochizuki-iut3}: simply because our two constructions  each provide a set  (denoted here by $\ttheta$) as a subset of  different ambient sets.
\item All that is claimed here (in Theorem~\ref{th:main}) is that there is a set $\ttheta$  of similar sort (which is constructed here quite explicitly) and that there is a lower bound  of a similar sort for its size. 
\eenum

\para Let me now come to another aspect of Arithmetic Teichmuller Theory and its relationship to \iut. An important technical device which plays a central role in \iut, is the notion of \emph{$\log$-links}. 

In \ssep\ref{se:log-links}, I prove the existence of a version of  Mochizuki's $\log$-links  from the point of view of  Arithmetic Teichmuller theory of \cite{joshi-tempered}. Notably (in Theorem~\ref{th:log-link}) I demonstrate that $\log$-links arise in \cite{joshi-tempered} naturally via \cite{fargues-fontaine} and that a $\log$-link  explicitly links two possibly topologically inequivalent algebraically closed perfectoid fields (see  Corollary~\ref{cor:log-link-in-arithmetic}). 

An important consequence of this is that the idea of moving between different untilts of $\cpt$ considered in \cite{joshi-tempered} is essentially equivalent to Mochizuki's idea of moving between different ``arithmetic holomorphic sturctures'' via log-links considered \iut. One can say that from the point of view of \cite{joshi-tempered}, that an untilt of $\cpt$ provides an (arithmetic) Berkovich analytic (=holomorphic) structure in the same sense that the codomain of a log-link of \iut\ provides a new arithmetic holomorphic structure in \iut.

\para The present document is almost self contained except for \cite{joshi-tempered} and neither use any results of \iut.

\para {\bf Acknowledgments} 
I would like to thank Matthew Morrow for comments on \ssep\ref{se:ansatz}. I am also indebted to Taylor Dupuy for reading and providing some comments on early drafts of this manuscript and for  conversations regarding \iut. I would also like to thank Kiran Kedlaya for some correspondence in the the early stage of this work and for his encouragement. Another mathematician who helped me during the course of the work presented here has chosen to remain anonymous, but I would like to acknowledge our correspondence during this the course of this work here.  As I have already stated in \cite{joshi-tempered}, neither that paper nor this one, could have existed without Mochizuki's work in anabelian geometry, $p$-adic Teichmuller Theory and his work on the Arithmetic Szpiro Conjecture. So my intellectual debt to Shinichi Mochizuki should be obvious. 

\section{Value groups of untilts}

\para\label{pa:value-group-comp-loc} Let me begin by noting that elements of value groups of all untilts of a fixed perfectoid field $F=\cpt$ can be compared in the value group of the tilt $F$.
\bpro 
Let $K_1$, $K_2$ be untilts of $\cpt$. . 
\benumlab
\item Then one has $\abs{K_1^*}=\abs{{\C_p^{\flat*}}}$. 
\item Hence, as one has $K_2^\flat=\cpt$,  so $K_1$ and $K_2$ have the same value group.
\item In particular, one can compare elements of $\abs{K_1^*}$ and $\abs{K_2^*}$ as elements of $\abs{{\C_p^{\flat*}}}$. 
\eenum
\epro 
\bp 
The first assertion is  \cite[Lemma 1.3.3]{kedlaya15} or \cite[Lemma 3.4(iii)]{scholze}: if $K_1^\flat=\cpt$ then one has $$\abs{K_1^*}=\abs{K_1^{\flat*}}=\abs{{\C_p^{\flat*}}}.$$ 
So if $K_1^\flat=K_2^\flat$ then $\abs{K_1^*}=\abs{K_2^*}=\abs{{\cpt}^*}$. So one can compare elements of value groups of $K_1$ and $K_2$ as elements of the value group of $\cpt$. This proves the proposition.
\ep

\section{Fargues-Fontaine curves $\syQp$ and $\sxQp$ and the ring $B$}
\para I will write $F=\cpt$, $\O_F\subset F$ be its valuation ring and $\fm_F\subset \O_F$ be its maximal ideal. This will be fixed throughout this note. Let $\sG/\Z_p$ be a Lubin-Tate formal group with logarithm $\sum_{n=0}^\infty \frac{T^{p^n}}{p^n}$. Then $\sG(\O_F)$ is naturally a Banach space over $\Q_p$ \cite[Chap. 4]{fargues-fontaine}.   Write $B=B_{\Q_p}$ for the ring constructed in \cite[Chap 2]{fargues-fontaine} for the datum $(F=\cpt,E=\Q_p)$.

\para Let $\sxQp$ be the complete Fargues-Fontaine curve, and  $\syQp$ be the incomplete Fargues-Fontaine curve \cite[Chapter 5]{fargues-fontaine}. Of interest to us are the residue fields of closed points $x\in\sxQp$ (resp. $y\in\syQp$). One has a canonical bijection of closed points of degree one \cite{fargues-fontaine}:
$$\abs{\syQp}/\vphi^\Z\to \abs{\sxQp}.$$
If $y\in \abs{\syQp}$ maps to $x\in\sxQp$ under this bijection then one has a natural isomorphism of residue fields
$$K_y\mapright{\isom} K_x.$$

By \cite[Section 10.1]{fargues-fontaine} one has action of 
$\gal(\bQ_p/\Q_p)$ (and hence an action of any of its open subgroups) on $\abs{\syQp}$ (resp. on $\sxQp$). 

By \cite[Theorem 8.29.1]{joshi-tempered} one  also has a natural an action of ${\rm Aut}_{\Z_p}(\sG(\O_F))$  (resp. ${\rm Aut}_{\Q_p}(\sG(\O_F))$) arising from the natural identification ( \cite[Section 2.3]{fargues-fontaine} (resp. \cite[Th\'eor\`eme 5.2.7]{fargues-fontaine})) of 
$$\abs{\syQp} \isom \left(\sG(\O_F)-\{0\}\right)/\Z_p^*\qquad \text{ resp. } \qquad \abs{\sxQp}\isom \left(\sG(\O_F)-\{0\}\right)/\Q_p^*.$$

\para 
For simplicity I will write $V=\sG(\O_F)$ considered as a topological $\Z_p$-module  and also as a $\Q_p$-Banach space. Both are equipped with an action of $\gal(\bQ_p/\Q_p)$ and in this notation one has 
$$\abs{\sxQp}\isom \P(V), \text{ and}$$
$$\abs{\syQp}\isom (V-\{0\})/\Z_p^*.$$

In this notation the automorphism groups are $${\rm Aut}_{\Z_p}(\sG(\O_F))={\rm Aut}_{\Z_p}(V)$$ and $${\rm Aut}_{\Q_p}(\sG(\O_F))={\rm GL}(V).$$

\para 
The $\gal(\bQ_p/\Q_p)$ action on $\sxQp$ has one unique fixed point   with (resp. with a finite orbit for any open subgroup) \cite[Proposition 10.1.1]{fargues-fontaine}, denoted by $\infty\in\abs{\sxQp}$ whose residue field is identified with $\C_p$ (with the standard valuation $v_{\C_p}(p)=1$). Similarly one has a canonical point $(\infty,\vphi=1)\in\abs{\syQp}$ whose image in $\sxQp$ is $\infty\in\sxQp$.

So one has two distinct actions 
$$\gal(\bQ_p/\Q_p)\act \abs{\sxQp} \curvearrowleft {\rm Aut}_{\Q_p}(\sG(\O_F))$$
which move points on the curve and hence their residue fields.

\para As was noted in \cite[\ssep 9.3]{joshi-tempered}, the theory of \iut\ is naturally a multiplicative theory. In the multiplicative context of \iut\ one can use the isomorphism
$$\sG(\O_F)\isom \hgm(\O_F).$$

Notably one has $${\rm Aut}_{\Q_p}(\sG(\O_F)) \isom {\rm Aut}_{\Q_p}(\hgm(\O_F))$$
and especially one has the multiplicative action:
$$\gal(\bQ_p/\Q_p)\act \abs{\sxQp} \curvearrowleft {\rm Aut}_{\Q_p}(\hgm(\O_F)).$$

The similarities and differences between this and \iut\ is discussed in \cite[Section 9]{joshi-tempered}.
\section{The Tate curve set up}
\para  Let $E$ be a $p$-adic field , $G_E$ be its absolute Galois group for a given algebraic closure of $E$.

\para Fix a geometrically connected, smooth quasi-projective variety $X/E$ where $E$ is a $p$-adic field. Then the objects of the Arithmetic Teichmuller space $\fjxe$ (resp. $\fjxe_{\cpt}$) constructed in \cite{joshi-teich} are triples $(Y/E',E'\into K)$ where $E'$ is a $p$-adic field, $Y$ is a geometrically connected, smooth quasi-projective variety over $E'$ with an isomorphism $\pit{Y/E'}\isom \pit{X/E}$ and $K$ is an algebraically closed perfectoid field containing $E'$ (resp. algebraically closed perfectoid field $K$ containing $E'$ and $K^\flat=\cpt$). The properties of $\fjxe$ (resp. $\fjxe_{\cpt}$) are described in \cite[\ssep 1.4]{joshi-tempered}.

\para The case of interest here is this: let $C/E$ be  an elliptic curve with split multiplicative reduction over $E$. Write $X=C-\{O\}$ be the standard elliptic cyclops. Let $q=q_C\in E^*$ be its Tate parameter; later on one may want to assume that $X/E$ is a a once punctured elliptic curve i.e. an elliptic cyclops which is defined over a number field i.e. $X/E$ is of strict Belyi type in the sense of \topics. In that case, as was noted in \cite{joshi-tempered}, for any $(Y/E',E'\into K)$ one has $Y\isom X$ as $\Z$-schemes.

\para Tate's Theory provides the following:
\bpro Let $\Delta_C$ be the minimal discriminant of $C/E$ and let $q$ be the Tate parameter.
One has 
$$\Delta_C=\Delta(q)=q\prod_{n=1}^\infty (1-q^n)^{24}.$$
\epro
\bp 
This is \cite[Chapter V, Theorem 3.1]{silverman-advanced}.
\ep

\para Assume from now on that the $p$-adic field $E$ contains $\Q_p(q^{1/2\ell},\zeta_{2\ell})$ for a chosen prime $\ell\neq p$ and $\ell\geq 3$.

\section{The choice of a theta function}
\para I will  choose a theta function on $C/E$. For the theory of $p$-adic theta-functions see \cite{roquette-book}.  A theta-function on $C$ is a function on the the universal cover with some quasi-periodicity properties--equivalently a theta-function is a section of a some line bundle.
Let $q_E$ be the a Tate parameter for $C/E$.

The function I choose here, namely $\tm u$   is the same as the function used in \cite[Proposition 1.4]{mochizuki-theta} and \iut:
\be
\tm u = q_E^{-1/8}\sum_{n\in\Z}q_E^{\frac{1}{2}(n+\frac{1}{2})^2}u^{2n+1}.
\ee
This satisfies the following properties
\benumlab
\item $\tm{u^{-1}}=-\tm u$,
\item $\tm {q_E^\frac{j}{2}u}=(-1)^j q_E^{-\frac{j^2}{2}}u^{-2j}\tm u$.
\eenum
Since one is working with $\sqrt{q_E}$ as opposed to $q_E$, this theta function naturally lives on a double cover $ C'\to C$ of $C$.

From the first formula one sees that $$\tm 1 =0,$$ using this with the quasi-periodicity given by (2) one sees that
$$\tm {q_E^\frac{j}{2}}=0 \text{ for all } j\in \Z,$$
and that $\tm u=0$ if and only if $u=q_E^{j/2}, j\in\Z$.
From this one can also provide a multiplicative description of this theta function along the lines of \cite{roquette-book} or \cite{silverman-advanced}. Secondly the points $q_E^{j/2}$ $(j\in\Z)$ are mapped to the identity element $O\in C'$ under the Tate parametrization of $C'$. 

\para There is a line bundle $\sL_\theta$ on $C'$ corresponding to $\tm u$. This line bundle may be trivialized on $Y=C'-\{O\}$ and one sees from the above that  $\tm{u}$ has no zeros on $Y$ and can be used to trivialize $\sL_\theta$ on $Y$. I will colloquially speak of the $\theta$-function $\tm{u}$ as a function on $Y$ (and by abuse of terminology, sometimes even on $C'$). This ties up with \cite[\ssep 1.2]{joshi-tempered}.

\para  Let $K$ be an algebraically closed perfectoid field equipped with an embedding of $$\iota:E\into K$$ of valued fields and work over $C_K^{an}$. The theory of Tate elliptic curves also says that $C/K$ is also a Tate elliptic curve and the Tate parameter $q$ of $C/K$ is given by $q=\iota(q_E)$ where $q_E$ is the Tate parameter of $C/E$. In particular the $\theta_M$ can be viewed as a $\theta$-function over $C_K^{an}$ and if one needs to remember $K$ then one writes $\theta_{M;K}$.

\section{Choice of theta values}\label{se:theta-values}
\para Let $\ell\geq 3$ be a prime with $\ell\neq p$  and let $\zeta_\ell\in\Qpb$ be a primitive $\ell^{th}$-root of unity.  Let $q^{1/2\ell}\zeta_\ell^\alpha$,  for $\alpha=0,\ldots,\ell-1$,  be the $\ell^{th}$-roots of $\sqrt{q}$. Then the following collection of elements will be referred to as $\theta$-values (for the given data):
$$\frac{1}{\xi_j}=\frac{\tm {{q^\frac{j}{2\ell}\zeta_\ell}}}{\tm {\zeta_\ell}}=(-1)^j {q^{-\frac{j^2}{2\ell}}}\zeta_\ell^{-2j}.$$
The theta values of interest to us are the $\theta$ values $\xi_j$ for $j=1,\ldots,\ells$. There is an evident action of $G_E$ on the set of values considered above. These values are obviously defined in any  algebraically closed perfectoid field $K\supset E$.

By the description of the torsion points of a Tate curve in \cite{serre-abelian}, the images of $q^{1/2\ell}\zeta_\ell^\alpha$ under the Tate parameterization of $C'$ provide a subgroup of order $\ell$ in $C'[\ell]$ on the double cover of $C'\to C$. 

\para Note that  for any algebraically closed perfectoid field $K\supset E$ one has
$$\xi_j\subset\O_{\bE}\subset \O_K.$$ Note that strictly speaking I should write $\theta_{M;K}(q^{1/2\ell}\zeta^j_{2\ell})$ etc.,  to indicate the dependence on $K$ and strictly speaking the values of interest to us are the ones calculated in $K$--namely:
$$\frac{\theta_{M;K} ({{q^\frac{j}{2\ell}\zeta_\ell}})}{\theta_{M;K} ({\zeta_\ell})}.$$
 I hope that my notational simplification $\xi_j$ will not be too confusing.

\para Let $K$ be an algebraically closed perfectoid field equipped with an embedding $E\into K$. Then one can, by base extension from $C^{an}/E$ to $K$, consider the analytic space $C^{an}_K$. The $\theta$-function considered above can now be considered in $C^{an}_K$. As $C/K$ is a Tate elliptic curve, it is described by a theory of $p$-adic $\theta$-functions with values in $K$ (by \cite{roquette-book}). In particular one can consider $\theta_M$ as a $\theta$-function over $C^{an}_K$ and for any $q,u\in K$ one may evaluate the theta function $\theta_M$ at these values in $C^{an}_K$. By \cite[Theorem~3.15.1]{joshi-tempered} one knows that if $K_1,K_2$ are two algebraically closed perfectoid fields which are not topologically isomorphic then  the analytic function theories of $C_{K_1}^{an}$ and $C_{K_2}^{an}$ are not isomorphic (as the these two analytic spaces are not isomorphic).   This allows us to talk about values of the $\theta$-functions, even at elements  of $E$ (a finite extension of $\Q_p$), as living in two topologically distinct perfectoid fields. To put it in the parlance of \iut, two topologically non-isomorphic algebraically closed perfectoid fields $K_1,K_2$ provide, quite literally, two distinct ``arithmetic holomorphic (i.e Berkovich analytic) structures'' (namely $C_{K_1}^{an}$ and $C_{K_2}^{an}$). By \cite{joshi-tempered} the values may be lifted to the ring $B$ where one can compare them.

\section{Mochizuki's ansatz ``$q=q^{j^2}$'' and the $\theta$-link via $\sxQp$}\label{se:ansatz}
\para Let $q$ be the Tate parameter  of a Tate elliptic curve over a $p$-adic field. According to \cite[Page 4]{mochizuki-gaussian}, Mochizuki's theory of $\theta$-links  is roughly based on the idea that at some level in his theory \iut\ one can declare that   (I am putting this in quotes intentionally) $$\text{``}q=q^{j^2}\text{''}.$$ 

Let $\ell\geq 3$ be a prime and let $\ells=\frac{\ell-1}{2}$ and let $q$ be the Tate parameter of a given Tate elliptic curve. In \iut\ Mochizuki, in fact, asserts that one can set
$$\text{``}\left(q^{1^2},q^{2^2},q^{3^2},\cdots, q^{\ells^2}\right)=q\text{''}$$
in his theory! \emph{This is the $\theta$-link of Mochizuki's theory \iut}. 

It is obviously enough to  be able to set
$$\text{``}p=\left(p^{1^2},p^{2^2},p^{3^2},\cdots, p^{\ells^2}\right)\text{''}.$$ 

\para My point here is that this can be done explicitly in the Arithmetic Teichmuller Theory of \cite{joshi-tempered} as I will now demonstrate. 

Let $a\in\fm_F$ and assume $a\neq0$. Then for $j=1,\ldots,\ells$ one has primitive elements of degree one  given by $$[a^{j^2}]-p\in W(\O_F)$$ and  they generate prime ideals  (see \cite[Section 2.2]{fargues-fontaine})
$$\wp_j=([a^{j^2}]-p)\subset W(\O_F),$$ 	
and these define closed points $y_1,\ldots,y_{\ells}$ of $\sxQp$, with residue fields $K_j=K_{y_j}$.
From now on these letters and symbols will be used with this meaning.

\para Let \be\Sigma_F=\left\{ ([a]-p,[a^{2^2}]-p,[a^{3^2}]-p, \cdots ,[a^{\ells^2}]-p): a\in\fm_F-\{0\}\right\}.\ee 
By \cite[Lemma 2.2.14]{fargues-fontaine} for each $1\leq j\leq \ells$, the ideal $$\wp_j=([a^{j^2}]-p)\subset W(\O_F)$$ is a principal prime ideal of $W(\O_F)$ and its extension to $W(\O_F)\subset B$ provides a closed maximal ideal of $B$ , and hence each provides a closed point of $\sxQp$ of degree one .

So any element $([a]-p,[a^{2^2}]-p,[a^{3^2}]-p, \cdots ,[a^{\ells^2}]-p)\in\Sigma_F$ defines a tuple of closed points $(y_1,\ldots,y_\ells)\in \sxQp^\ells$ with tuple of residue fields $(K_j)_{1\leq j\leq\ells}=(K_{y_j})_{1\leq j\leq\ells}$. I will usually think of an element of $\Sigma_F$ as providing a tuple of points of $\sxQp$.  \textcolor{red}{Note that an arbitrary tuple in $\sxQp^{\ells}$ need not belong to $\Sigma_F$.} So the set $\Sigma_F$ should not be confused with $\sxQp^\ells$.

\para Notably, let $t\in \cpt$ be such that $([t]-p)\subset W(\O_{\cpt})$ gives the canonical point of $\sxQp$ (\cite[Chap 10]{fargues-fontaine}) with residue field $\C_p$ as a $\gal(\bQ_p/\Q_p)$-module. Then the tuple
\be\label{eq-canonical-tuple}
  ([t]-p,[t^{2^2}]-p,[t^{3^2}]-p, \cdots ,[t^{\ells^2}]-p)\in\Sigma_F
\ee
 in which the first factor is corresponds to the canonical point of $\sxQp$.
 
\para\label{pa:Mochizuki-Ansatz} When one works with a fixed perfectoid field $K\supset \Q_p$ with $K^\flat=\cpt$ it is conventional to write the principal ideal of $W(\O_F)$ defining $\O_K$ (and hence $K$) as $([p^\flat]-p)\subset W(\O_F)$ i.e. $([p^\flat]-p)=\ker(\eta_K:W(\O_F)\onto \O_K)$. In this notation one has a  point of $\Sigma_F$ given by the system of equations
\begin{align*}
[p^{\flat}]-p&=0,\\
[(p^{\flat})^{2^2}]-p&=0\\
\vdots &\quad\ \vdots\\
[(p^{\flat})^{\ells^2}]-p&=0.
\end{align*}
And hence one has
$$p=[p^{\flat}]=[(p^{\flat})^{2^2}]=
\cdots=[(p^{\flat})^{\ells^2}].
$$
Proposition~\ref{pr:lift-vals} given below shows how valuations of the corresponding tuple of perfectoid fields $(K=K_1,K_2,\ldots,K_\ells)$ are also simultaneously scaled by suitable factors--\emph{this simultaneous scaling of valuations is the primary reason why this ansatz  is  important in \iut--at least as far as I understand.}
 
\para There is a natural continuous action of $\gal(\bQ_p/\Q_p)$ on $\Sigma_F$ which agrees with the natural $\gal(\bQ_p/\Q_p)$ action on the set of primitive 
elements of degree one of $W(\O_F)$. On a tuple $(y_1,\ldots,y_{\ells})$ this action is given by
$$y_j\mapsto \sigma(y_j)=[\sigma(a^{j^2})]-p=\sigma([a^{j^2}]-p).$$
so that $(\sigma(y_1),\ldots,\sigma(y_{\ells}))\in\Sigma_F$ has  similar valuation theoretic properties as  $(y_1,\ldots,y_\ells)$. Hence one has proved the following assertion:

\blem 
The set $\Sigma_F$ is Galois stable, and notably the canonical tuple of \eqref{eq-canonical-tuple} has the property that the image of  its first coordinate in $\sxQp$ is the canonical point of $\sxQp$ and hence is fixed  under the action $\gal(\bQ_p/\Q_p)$.
\elem

\para\label{pa:norms-vs-element-comp} Note that as each of the fields $K_1,\ldots,K_\ells$ are algebraically closed perfectoid and hence  each field $K_j$ contains a copy of the algebraic closure of $E$,  $\bE_j\subset K_j$. But by \cite{kedlaya18} there may be no topological isomorphism $K_1\isom K_2$ (say) which takes $\iota_1:\bE_1\into K_2$ isomorphically to $\iota_2:\bE_2\into K_2$. Thus one must take care when working with elements of our field $E$ and its algebraic closures in the fields $K_j$. However  because of \ssep\ref{pa:value-group-comp-loc} one may be able to compare the absolute values $\abs{\iota_1(z)}_{K_1}$ and $\abs{\iota_2(z)}_{K_2}$ but not necessarily compare $\iota_1(z)$ and $\iota_2(z)$. I will simply write $\abs{z}_{K_1}$ instead of $\abs{\iota_1(z)}_{K_1}$. Hopefully this will not cause any confusion.

\para One can achieve  Mochizuki's ansatz $$\text{``}\left(q^{1^2},q^{2^2},q^{3^2},\cdots, q^{\ells^2}\right)=q\text{''}$$ and also achieve simultaneous scaling of valuations using  by working with elements of Mochizuki's Ansatz $\Sigma_F$ as I will now demonstrate. Here is how this takes place.
\bthm\label{pr:lift-vals} 
Let $(y_1,\ldots,y_\ells)$ be closed points of $\sxQp$ corresponding to an element of $\Sigma_F$. Then one has (remembering \ssep\ref{pa:norms-vs-element-comp})
\benumlab
\item $$v_{K_j}(p)=j^2v_{K_1}(p)\text{ for } j=1,\ldots,\ells.$$
\item Notably if $K_1=\C_p$ then $v_{K_j}(p)=j^2$.
\item Hence for  $z\in \bE\subset K_j$, for $j=1,\ldots,\ells$ one has
$$\abs{z}_{K_j}=\abs{z}_{K_1}^{j^2}.$$
\item However the valued fields $K_1,\ldots,K_\ells$ need not be all topologically isomorphic.
\eenum
\ethm

\bp 
\textcolor{red}{Note the since one can view the values $\abs{-}_{K_j}$ in the value group of $F$  by \ssep\ref{pa:value-group-comp-loc} one can make such assertions about $\abs{-}_{K_j}$.}

The valuation of $K_j$ can be computed from \cite[Proposition 2.2.17]{fargues-fontaine} as $$v_{K_j}(p)=v_F(a^{j^2})=j^2v_F(a)=j^2v_{K_1}(p).$$
Hence the first assertion is proved and the second follows from this. Note that by loc. cit. the restriction of $\abs{-}_{K_j}$ to $E\subset K_j$ is given by the above formula and hence the penultimate assertion is also established. The last claim is an immediate  consequence of \cite{kedlaya18}.
\ep

\brem 
Let me remark that I use the projection to the first coordinate $$\Sigma_F \subset \sxQp^\ells \to\sxQp$$ for valuation computations. There are also more complicated (and highly non-algebraic) maps from $\Sigma_F\to \sxQp$ possible for example $$([a]-p,[a^{2^2}]-p,[a^{3^2}]-p, \cdots ,[a^{\ells^2}]-p)\mapsto [a^{\sum_{j=1}^\ells j^2}]-p.$$
\erem

\brem
I do not expect that $\Sigma_F$ is stable under ${\rm Aut_{\Q_p}(\sG(\O_F))}$. It seems reasonable, following Mochizuki's ideas in \cite{mochizuki-iut3}, that one should work with ${\rm Aut}(\gal(\bQ_p/\Q_p)\act\sG(\O_{\cpt}))$ instead and consider $\widetilde{\Sigma}_F$ which contains $\Sigma_F$ and is stable under ${\rm Aut}(\gal(\bQ_p/\Q_p)\act\sG(\O_{\cpt}))$.
\erem

\brem 
The set $\Sigma_F$ satisfying these properties should be viewed as providing arithmetic Kodaira Spencer classes. From \cite[proof of Theorem 3.15.1]{joshi-tempered} of \cite{kedlaya18} one sees that these are non-trivial because as one moves over $\Sigma_F$ using the action of the above groups,  the fields $K_1,\ldots,K_\ells$ need not be topologically isomorphic in general.
\erem

\newcommand{\perf}{\mathcal{Perf}}
\brem
The following remark is useful to remember. Let $\perf_{\cpt}$ be the category of algebraically closed perfectoid fields of characteristic zero with tilt isometric to $\cpt$ in which morphisms are morphisms of valued fields. Then one has a functor $$\fjxe_{\cpt}\to \perf_{\cpt}$$ which maps $(Y/E',E'\into K)\mapsto K$. 
Now consider the product category $\perf_{\cpt}^\ells$. Then $\Sigma_F$ is naturally a subset of objects of $\perf_{\cpt}^\ells$ and one can consider the full subcategory of $\perf_{\cpt}^\ells$ whose objects are $\Sigma_{\cpt}$. I will habitually identify the set $\Sigma_F$ with this subcategory of $\perf_{\cpt}$. Now let $$\fjxe^{Ansatz}_{\cpt}\to \perf_{\cpt}^\ells$$ be the category whose objects are $(Y/E', K_1,\ldots,K_\ells)$ where $\pit{Y/E'}\isom \pit{X/E}$ and $(K_1,\ldots, K_\ells)$ is the tuple of algebraically closed perfectoid fields arising as the tuple of residue fields of  a point of $\Sigma_F\into \perf_{\cpt}^\ells$. Strictly speaking I should work with $$\fjxe^{Ansatz}_{\cpt}\to \Sigma_{\cpt}\into \perf_{\cpt}^\ells$$ instead of the set $\Sigma_F$. 
Notably the construction of the set $\ttheta$ which is outlined in the next few sections strictly speaking takes place over $\fjxe^{Ansatz}_{\cpt}\to \Sigma_{\cpt}$ but I will simply conflate this with the set $\Sigma_{\cpt}=\Sigma_F$. Hopefully readers will have no trouble making the translation.
\erem

\section{Lifting theta values to $B$--Existence of Teichmuller lifts}\label{pa:lifting-values-to-B}
\para To use above ansatz, and work with the $\theta$-values (\ssep\ref{se:theta-values}) uniformly, it is best to work with lifts to the Fargues-Fontaine ring $B$. This allows us to exploit the various group actions on $\sxQp$ explicitly.

\para To understand the lifts of these values to $B$, I will use the following proposition. For a perfectoid field $K$ with $K^\flat=\cpt=F$ write $\eta_K:W(\O_F)\to \O_{K}$ for the canonical surjection \cite{fontaine94a} or \cite{fargues-fontaine}. This surjection is traditionally denoted by $\theta$ so the change of notation is naturally forced upon us. For $\rho\in(0,1)\subset \R$ let $\abs{-}_\rho$ be the multiplicative norms defining the Fr\'echet structure of $B$. Let $T_K\subset \bpip\subset B$ be the Tate module of $T_K\subset \sG(\O_K)\isom \bpip$. For a closed point $y$ of degree one of $\sxQp$ with $K=K_y$, I will write $T_y=T_{K_y}$ for simplicity. \emph{Note that $T_y$ is a free $\Z_p$-module of rank one.}

\bpro\label{pr:teichmuller-lift-A}
Let $\xi$ be one of the $\theta$-values. Let $y$ be a closed point of degree one of $\sxQp$ with residue field $K_y$. Then 
\benumlab 
\item there exists $x\in\fm_F\subset \O_F$ such that 
$\eta_{K_y}([x])=\xi$.
\item $\abs{[x]}_\rho=\abs{x}_F=\abs{\xi}_{K_y}$,
\item and for any $\tau\in T_y\subset \bpip$ one has 
$$\eta_{K_y}(\tau+[x])=\xi.$$
\item So the values $\tau+[x]$ provide a lift of $\theta$-values $\xi$ and notably one has
$$\abs{[x]}_\rho\leq \sup\left\{\abs{\tau+[x]}_\rho:\tau\in T_y\right\}.$$
\eenum
\epro
\bp 
The first assertion has nothing to do with specifics of theta values $\xi$ and is a consequence of \cite[Corollary 2.2.8]{fargues-fontaine} which asserts that there exists a Teichmuller lift of every element of  $ \O_{K}$ under the canonical surjection $\eta_K:W(\O_F)\to\O_K$. The second assertion uses just the definition of $\abs{-}_{K_y}$ given by \cite[Proposition 2.2.17]{fargues-fontaine} and (1). The third property is immediate from the fact that $\ker(\eta_{K_y})\supset T_y$, in fact, one has the exact sequence of Banach spaces (\cite[Propostion 4.5.14]{fargues-fontaine})
$$\xymatrix{0\ar[r] & T_y\tensor \Q_p\ar[r] &  \bpip \ar[r]^{\eta_{K_y}}& K_y\ar[r] & 0.}$$
The fourth property is self-evident.
\ep

\brem 
Let $[x]$ be a  lift of $\xi$ to $B$ with $\eta_{K_y}([x])=\xi$. Then one is interested in bounds for $\abs{[x]}_\rho=\abs{\xi}_{K_y}$. The idea is to instead consider bounds for elements of $T_y+[x]$.

Note that union of all such lifts $T_y+[x]$ of $\xi$ is not closed under Frobenius. So one needs to enlarge the locus to be Frobenius stable and then consider its intersection with $\bpip$.
\erem

\bpro\label{pr:teichmuller-lift-B}
Let $(y_1,\cdots,y_\ells)$ be a tuple of closed points of $\sxQp$ corresponding to an element of $\Sigma_F$. Let $[x_j]\in W(\O_F)\subset B$ be a Teichmuller lift of $\xi_1$  under
$\eta_{K_j}:W(\O_F)\to \O_{K_j}$ then for $j=1,\ldots,\ells$. Then one has
$$\abs{[x_j]}_\rho=\abs{[x_1]}^{j^2}_\rho.$$
\epro
\bp 
This is clear from Proposition~\ref{pr:lift-vals} and the fact that $[x_j]$ is a lift of $\xi_1$ under $\eta_{K_j}$, and hence $$\abs{[x_j]}_\rho=\abs{\xi_1}_{K_j}=\abs{\xi_1}_{K_1}^{j^2}=\abs{[x_1]}_\rho^{j^2}.$$
\ep

\section{Construction of $\ttheta\subset B^\ells$}\label{se:construction-ttheta}
\para Now one can start the construction of a subset of interest. Let me begin with the following consequence of the results of the previous section.

\bthm\label{thm:theta-pilot-object-appears} 
Let $\xi=\xi_1$ be the theta value constructed in \ssep\ref{se:theta-values}. Let $(y_1,\cdots,y_\ells)$ be a tuple of closed points of $\sxQp$ corresponding to an element of $\Sigma_F$. Let $(K_1,\ldots,K_{\ells})$ be the tuple of (perfectoid) residue fields of $(y_1,\cdots,y_\ells)$. Let $[x_j]$ be a Teichmuller lift  to $B$ of $\xi_1$ under $\eta_{K_j}$. Then
\benumlab
\item
\be
 \prod_{j=1}^{\ells}\abs{[x_j]}_\rho=\prod_{j=1}^{\ells}\abs{\xi}_{K_1}^{j^2}, 
\ee
\item 
or additively this is
\be \label{eq:log-theta-pilot}
\sum_{j=1}^{\ells}\log\abs{[x_j]}_\rho=\sum_{j=1}^{\ells}j^2\cdot\log\abs{\xi}_{K_1}.
\ee
\item 
Notably in the left hand side of \eqref{eq:log-theta-pilot},  one can replace the Teichmuller lifts  $\abs{[x_j]}_\rho$ by $\abs{[x_j]+\tau_{j}}_\rho$ where $\tau_j\in T_{y_j}$ in the above quantity and work with the supremum of all such values:
\be 
\sup_{\tau_1,\ldots,\tau_{\ells}}\left\{\sum_{j=1}^{\ells}\log\abs{[x_j]+\tau_j}_\rho\right\} \geq \sum_{j=1}^{\ells}\log\abs{[x_j]}_\rho=\sum_{j=1}^{\ells}j^2\cdot\log\abs{\xi}_{K_1}.
\ee
where the supremum is taken over $\tau_j\in T_{y_j}$ for $j=1,\ldots,\ells$.
\item Notably for $K_1=\C_p$ one has
\be\label{eq:eq:log-theta-pilot-variation} 
\sup_{\tau_1,\ldots,\tau_{\ells}}\left\{\sum_{j=1}^{\ells}\log\abs{[x_j]+\tau_j}_\rho\right\} \geq \sum_{j=1}^{\ells}j^2\cdot\log\abs{\xi}_{\C_p}.
\ee
\eenum
\ethm

\bp 
The proof follows by putting together Proposition~\ref{pr:teichmuller-lift-A} and Proposition~\ref{pr:teichmuller-lift-B}.
\ep

\brem 
Let me remark that one has from the definition of $\xi_1,\ldots,\xi_{\ells}$ that
$$\abs{\xi_j}_{\C_p}=\abs{\xi_1}^{j^2}_{\C_p}$$ and hence one can write 
$$\sum_{j=1}^{\ells}\log\abs{\xi_j}_{\C_p}=\sum_{j=1}^{\ells}\log\abs{\xi_1}^{j^2}_{\C_p}=\sum_{j=1}^{\ells}{j^2}\cdot \log\abs{\xi_1}_{\C_p}.$$
\erem

\brem 
Let me put Theorem~\ref{thm:theta-pilot-object-appears} in the perspective of  \cite[Corollary 3.12]{mochizuki-iut3}. 
\benumlab
\item The   quantity 
\be\label{eq:approx-theta-pilot-object}
\frac{1}{\ells}\sum_{j=1}^{\ells}j^2\cdot\log\abs{\xi}_{\C_p}
\ee
appears in the definition of  the ``arithmetic degree of the hull of  $\theta$-pilot object''  in \cite[Corollary 3.12]{mochizuki-iut3} and that degree is a key ingredient in that corollary.
\item Mochizuki appears to be the first to recognize that this quantity \eqref{eq:approx-theta-pilot-object} can be subjected to a meaningful variation  arising from the additional freedom available from his anabelian reconstruction point of view.   In his work this variation is achieved by means of the \emph{``Yoga of Indeterminacies''}  and the assertion of \cite[Corollary 3.12]{mochizuki-iut3} is that the resulting quantity bounds from above $\abs{q^{1/2\ell}}_{\C_p}$.
\item In Theorem~\ref{thm:theta-pilot-object-appears} the element of $\Sigma_F$ giving rise to $K_1,\ldots,,K_\ells$ can be considered to be a variable in  \eqref{eq:log-theta-pilot}. Thus it makes perfect sense to vary the element of $\Sigma_F$  and consider all the values so obtained. The supremum of these values will bound any specific value such as $\sum_{j=1}^{\ells}j^2\cdot\log\abs{\xi}_{\C_p}$. 
\item One can also enlarge the set of lifts $[x_j]+\tau_j\in B$ by replacing it by a Frobenius and Galois stable subset of $B$ generated by the $[x_j]+\tau_j$.
\item \emph{Strictly speaking, for Diophantine applications, one should do this for all primes of semi-stable reduction simultaneously.} So this should be considered to be a prototype for global constructions.
\item  Note that the terms on the left hand side of \eqref{eq:log-theta-pilot} are evidently mixed by the actions of  symmetries available in Arithmetic Teichmuller Theory of \cite{joshi-tempered}.
\item Theorem~\ref{thm:theta-pilot-object-appears} shows that a similar variation of \eqref{eq:approx-theta-pilot-object} exists from the point of view of Arithmetic Teichmuller Theory of \cite{joshi-tempered}.  Note however that Mochizuki does not work with the ring $B$, and hence the set constructed here is not the same as the one used in \iut. \emph{So I do not claim that what is proved here also provides a proof of \cite[Corollary 3.12]{mochizuki-iut3}.}
However what is proved here does provide an independent proof of the important assertion of \iut\  that the  quantity \eqref{eq:approx-theta-pilot-object} is subject to a natural variation, achieved in the present theory via the Teichmuller Theory of \cite{joshi-teich} and the existence  of Mochizuki's Ansatz (\ssep\ref{se:ansatz}) which leads naturally to the left hand side of \eqref{eq:eq:log-theta-pilot-variation}, which may be used to bound \eqref{eq:approx-theta-pilot-object}.
\eenum 
\erem

\para\label{pa:bound-of-cor312} 
Let $\ttheta\subset (B)^\ells=B\times B\times \cdots\times B$ ($\ells$ factors) be the smallest possible subset
containing all the lifts of $\xi_1$ of the form $[x_j]+\tau_j\in B$ obtained from $\Sigma_F$ such that $\ttheta$ is both
 \benumlab
 \item is Galois stable, 
 \item and also ${\rm Aut}(G_E\act \sG(\O_F))$ stable. 
 \eenum
Then $\ttheta$ should be considered  to be substitute, in Arithmetic Teichmuller Theory of \cite{joshi-tempered}, for the hull of $\theta$-pilot object(s) in \cite[Corollary 3.12]{mochizuki-iut3}. A possible reformulation of that Corollary from my point of view  is this: 
\be\label{eq:inequality-of-cor3.12}|\ttheta|_B:=\sup_{\rho\in(0,1)}\left\{|\ttheta|_\rho\right\}=\sup_{\rho\in(0,1)}\sup\left\{\prod_{j=1}^\ells\abs{z_i}_\rho: z=(z_1,\cdots,z_\ells)\in \ttheta \right\} \geq \abs{q^{1/2\ell}}^\ells_{\C_p}.\ee
For other treatments of  \cite[Corollary 3.12]{mochizuki-iut3} from Mochizuki's point of view and how it may be used to prove a number of Diophantine conjectures  see \iut, or  \cite{fucheng}, \cite{yamashita}, \cite{dupuy2020statement} (where \cite[Corollary 3.12]{mochizuki-iut3} is stated as Conjecture 1.0.1), \cite{dupuy2020probabilistic} and \cite{dupuy2021kummer}.

\para Let me remark that the inequality \eqref{eq:inequality-of-cor3.12} is a version (in the theory of \cite{joshi-tempered}) of the inequality asserted in \cite[Corollary 3.12]{mochizuki-iut3}. To see that this is indeed the case note that if $0<x<1$ is a real number then $\log(x)<0$ and so by definition
 $$\abs{\log(x)}=-\log(x).$$
As $q$ is the Tate parameter one has $\abs{q}_{\C_p}<1$ and hence $\abs{q^{1/2\ell}}_{\C_p}<1$. Now suppose that $|\ttheta|_B<1$. For the bound of \eqref{eq:inequality-of-cor3.12} to be useful \emph{at all} this would certainly need to be the case--also see Proposition~\ref{pr:trivial-upper-bound-on-ttheta}. Then $$|\log(|\ttheta_B|)|=-\log(|\ttheta_B|)$$ and similarly $$|\log(|q^{1/2\ell}|_{\C_p}^\ells)|=-\log(|q^{1/2\ell}|_{\C_p}^\ells)$$
and so the inequality \eqref{eq:inequality-of-cor3.12} can also be written as
\be\label{eq:inequality-of-cor3.12-2}
-|\log(|\ttheta|_B)|\geq -|\log(|q^{1/2\ell}|_{\C_p}^\ells)|.
\ee
which is precisely the sort of inequality of \cite[Corollary 3.12]{mochizuki-iut3} especially if one writes $$|\ttheta_{Mochizuki}|_B=\sqrt[\ells]{|\ttheta|_B},$$ or equivalently
$$\log(|\ttheta_{Mochizuki}|_B)=\frac{1}{\ells}\log(|\ttheta|_B),$$
 then one gets from \eqref{eq:inequality-of-cor3.12-2} that
\be 
-|\log(|\ttheta_{Mochizuki}|_B)|\geq -|\log(|q^{1/2\ell}|_{\C_p})|.
\ee
which has the same shape as the inequality of \cite[Corollary 3.12]{mochizuki-iut3}. Since  Mochizuki works with a different set from the set $\ttheta$ constructed here so one cannot assert that it is exactly the same inequality. But the inequality stated here is the analog of Mochizuki's inequality in the theory of \cite{joshi-tempered}.

\section{The bound $|\ttheta|_B\geq |q^{1/2\ell}|^\ells_{\C_p}$}
\para 
Let me now prove how one can use $\Sigma_F$ to prove the lower bound:

\bthm\label{th:main} 
If $\ell$ is sufficiently large then one has $$|\ttheta|_B\geq |q^{1/2\ell}|^\ells_{\C_p}.$$
\ethm
\bp 
Let $t\in \cpt=F$ be such that $([t]-p)\subset W(\O_F)$ is the prime ideal corresponding to the canonical point of $\sxQp$. This means the residue field of this point is $\C_p$ with its natural action of $\gal(\bQ_p/\Q_p)$ and the valuation of $\C_p$ has the standard normalization $v_{\C_p}(p)=1$. 

Let $\ell\geq 3$ be a prime, $\ells=\frac{\ell-1}{2}$. Choose a root $t^{1/\ells^2}\in\cpt$. Let  for  $j=1,\cdots,\ells$,
$$a_j=\left(t^{1/\ells^2}\right)^{j^2}\in\cpt,$$
and prime ideals
$$\wp_j=[a_j]-p.$$
Then one has $$([a_1]-p,[a_2]-p,\ldots,[a_\ells]-p)\in\Sigma_F,$$
and let $K_1,\ldots,K_\ells$ be their residue fields and by construction $K_\ells$ is the canonical point of $\sxQp$ notably the residue field $K_\ells=\C_p$ (as a $\gal(\bQ_p/\Q_p)$-module). 

Now let us compute the relationship between valuations of $K_1,\ldots,K_\ells$. This can be done by Proposition~\ref{pr:lift-vals}. One has for $j=1,\ldots,\ells$:
$$v_{K_j}(p)=j^2\cdot v_{K_1}(p),$$
and hence in particular one has
$$v_{K_\ells}(p)=v_{\C_p}(p)=1=\ells^2 \cdot v_{K_1}(p),$$
and hence one has
$$v_{K_1}(p)=\frac{1}{\ells^2}\cdot v_{\C_p}(p).$$
So one has
\be\label{eq:val-rel-main} v_{K_j}(p)= j^2\cdot v_{K_1}(p)=\frac{j^2}{\ells^2}\cdot v_{\C_p}(p).
\ee

Now let $[x_j]\in W(\O_F)\subset B$ be a Teichmuller lift of $\xi=\xi_1\in K_j$ with $\abs{\xi}_{K_j}=\abs{x_j}_F$ (as in Proposition~\ref{pr:teichmuller-lift-A}, \ref{pr:teichmuller-lift-B}) and Theorem~\ref{thm:theta-pilot-object-appears}.
Then
\be\label{eq:val-rel-main2} 
\sum_{j=1}^\ells\log\abs{[x_j]}_\rho=\sum_{j=1}^\ells \log\abs{\xi}_{K_j},
\ee
and by the established relationship between valuations \eqref{eq:val-rel-main} one has
$$\abs{\xi}_{K_j}=\abs{\xi}_{\C_p}^{\frac{j^2}{\ells^2}}.$$
So one can write the right hand side of \eqref{eq:val-rel-main2} as
$$\sum_{j=1}^\ells\log\abs{[x_j]}_\rho=\sum_{j=1}^\ells \frac{j^2}{\ells^2}\log\abs{\xi}_{\C_p}.$$

This gives 

$$\sum_{j=1}^\ells\log\abs{[x_j]}_\rho=\sum_{j=1}^\ells \frac{j^2}{\ells^2}\log\abs{\xi}_{\C_p}=\frac{\ells(\ells+1)(2\ells+1)}{6\ells^2}\log\abs{\xi}_{\C_p}.$$

Now using $$\abs{\xi}_{\C_p}=\abs{q^{1/2\ell}}_{\C_p},$$ the right hand side of the above equation simplifies to
$$\frac{\ells(\ells+1)(2\ells+1)}{6\ells^2}\cdot\frac{1}{2\ell}\log\abs{q}_{\C_p}.$$
Using $\ells=\frac{\ell-1}{2}$, or equivalently $2\ells+1=\ell$ and canceling factors  this simplifies to
$$\frac{1}{12}\left(1+\frac{1}{\ells}\right)\log(\abs{q}_{\C_p}).$$
If $\ell$ is sufficiently large then $$\frac{1}{12}\left(1+\frac{1}{\ells}\right)<\frac{1}{4}\left(1-\frac{1}{\ell}\right)=\frac{\ells}{2\ell},$$ 
and as $\log(\abs{q}_{\C_p})<0$ one gets
$$\frac{1}{12}\left(1+\frac{1}{\ells}\right) \cdot \log(\abs{q}_{\C_p})> \frac{1}{4}\left(1-\frac{1}{\ell}\right)\cdot \log(\abs{q}_{\C_p})=\log(\abs{q^{1/2\ell}}^\ells_{\C_p}).$$
Thus $\Sigma_F$ contains a point at which $$\sum_{j=1}^\ells\log(\abs{[x_j]}_{\rho})>\log(\abs{q^{1/2\ell}}^\ells_{\C_p}).$$
Since $\log(|\ttheta|_B)$ is the supremum over $\Sigma_F$ (or over a bigger set), one obtains the asserted bound
$$|\ttheta|_B \geq \abs{q^{1/2\ell}}^\ells_{\C_p}.$$ 
In terms of $\ttheta_{Mochizuki}$ this is
$$\log|\ttheta_{Mochizuki}|\geq \log(|q^{1/2\ell}|_{\C_p}).$$
\ep

\bcor 
In the above notation, if $c\in \R^*$ is a non-zero real number such that $$|\ttheta|_B\leq c|q^{1/2\ell}|_{\C_p}^{\ells}$$
then $$c\geq 1.$$
\ecor 
\bp 
 If $c<1$ then this contradicts the lower bound provided by Theorem~\ref{th:main}.
\ep

\brem
This is essentially the version of \cite[Corollary 3.12]{mochizuki-iut3} valid in the present theory.
\erem

\brem 
It is possible that the set $\ttheta$ constructed above is too big--for obtaining non-trivial upper bounds. It is certainly enough to work with the (finite) union of
 $G_{E'}$-orbit of the point of $\Sigma_F$ constructed in the proof of Theorem~\ref{th:main} where $E'$ runs over all the finitely many extensions $E'\subset \bQ_p$ which are anabelomorphic with $E$ i.e. $G_{E'}\isom G_{E}$ (in the context of \iut, this corresponds to ``Mochizuki's Indeterminacy Ind1'') and replace $\ttheta$ by the smaller subset corresponding to this finite union of orbits.
\erem

\para The following proposition shows that it is not unreasonable to expect that $|\ttheta|_B$ is bounded from above. Let $\rho\in (0,1]\subset \R$ and let 
\be 
|\theta|_{B,\rho}=\sup_{z}\left\{ \prod_{j=1}^\ells \abs{z_i}_\rho: z=(z_1,\ldots,z_\ells)\in \ttheta \right\}.
\ee 
With this notation one has a trivial upper bound:
\bpro\label{pr:trivial-upper-bound-on-ttheta} 
For the set $\ttheta$ constructed in \ssep\ref{se:construction-ttheta} and for $\rho=1$, one has:
\be 
|\ttheta|_{B,1}\leq 1.
\ee
\epro
\bp 
Let $z=(z_1,\ldots,z_\ells)\in\ttheta$ and let $(y_1,\ldots,y_\ells)$ be the tuple of closed points of $\sxQp$ with the property that $\eta_{K_j}(z_j)=\xi_1$ i.e. $z_j$ is the lift of $\xi_1\in K_j$ to $B$ (here $K_j$ is the residue field of $y_j$). Let $\rho=1$ and consider the norm $\abs{-}_\rho=\abs{-}_1$ on $B$. 

Then I claim that for each $1\leq i\leq \ells$ one has $$\abs{z_i}_1\leq 1.$$
Indeed, any $z_i$  considered here is of the form $z_i=[x]+\lambda_i$ for some $\lambda_i\in T_{y_i}\subset B^{\vphi=p}$.  By \cite[Proposition 4.1.3]{fargues-fontaine} one has $$B^{\vphi=p}=(B^+)^{\vphi=p}\subset B^+\subset B$$ and by \cite[Proposition 1.10.7]{fargues-fontaine} one has $$B^+=  \left\{\lambda\in B: \abs{\lambda}_1\leq 1\right\}.$$ 
So $\lambda_i$ satisfies $\abs{\lambda_i}_1\leq 1$ and hence one sees that 
$$\abs{z_i}_1=\abs{[x]+\lambda_i}_1\leq \max\left(\abs{[x]}_1,\abs{\lambda_i}_1\right)\leq \max\left(\abs{[x]}_1,1\right).$$
But by the definition of $\abs{-}_1$ (\cite[D\'efinition 1.4.1]{fargues-fontaine}) one has
 $$\abs{[x]}_1 = \abs{x}_F <1 \text{ for any }x\in\fm_F.$$ 
 So $\abs{z_i}_1\leq 1$. This proves the claim.
\ep

\section{Mochizuki's $\log$-links}\label{se:log-links}
\para I will now demonstrate that a version Mochizuki's  $\log$-links are present in the theory of \cite{joshi-tempered}. Mochizuki's theory of $\log$-links is detailed in \iut--especially \cite[Page 24]{mochizuki-iut3} (and \topics); other treatments are available--\cite[Section 6.1, Definition 6.3]{fucheng}, \cite{yamashita} and \cite{dupuy2020statement}. I will not recall Mochizuki's Theory of $\log$-links as the approach provided below provides a cleaner demonstration of an  assertion crucial  to \iut\ namely that a  $\log$-link intertwines two possibly distinct (topological) field structures.

\para Consider $\hgm$ as a formal (Lubin-Tate) group over $\Z_p$ and let $E=\Q_p$, $F=\cpt$ and let $\sG/\Z_p$ be the Lubin-Tate group with formal logarithm $\sum_n\frac{T^{p^n}}{p^n}$. Let $\sgt$ be the universal cover of $\sG$ given by $\sgt$ as the limit of diagrams (see \cite[Chap. 4]{fargues-fontaine}) $$\sgt=\cdots\mapright{p}\sG\mapright{p} \sG\mapright{p}\sG.$$
For $\sG=\hgm$ this can be described quite concretely as follows: $\hgm(\O_{K})=1+\fm_K$ and so 
$$\widetilde{\hgm(\O_K)}=\left\{(x_n): x_n\in 1+\fm_K,  x_{n+1}^p=x_n  \right\}\isom \hgm(\O_{\cpt})=1+\fm_{\cpt}.$$
and for simplicity of notation, I will habitually write $\widetilde{1+\fm_K}=\widetilde{\hgm(\O_K)}\isom \hgm(\O_{\cpt})\isom 1+\fm_{\cpt}$. The mapping $$\hgm(\O_{\cpt})\to K$$ is the $p$-adic logarithm $$\hgm(\O_{\cpt})=1+\fm_{\cpt} \mapright{\log} K$$  given by $$1+\fm_{\cpt} \ni x\mapsto \log([x]),$$
where $[-]:\O_{\cpt}\to W(\O_{\cpt})$ is the Teichmuller lift (\cite[Example 4.4.7]{fargues-fontaine}). These constructions provide the following:
\bthm\label{th:log-link} 
Let $K_1,K_2$ be two characteristic zero untilts of $\cpt$.
There is a commutative diagram of continuous homomorphisms of topological vector spaces in which vertical arrows are isomorphisms
\be \xymatrix{
{} &	\sgt(\O_{K_1})\ar[r]^{\isom} \ar[d]^\isom &  \sG(\O_{\cpt})\ar [r]^{\log_\sG}\ar[d]^\isom & K_2\ar[d]^\isom \ar [r] & 0\\
\widetilde{1+\fm_{K_1}}\ar[r]^\isom &	\widetilde{\hgm(\O_{K_1})}\ar[r]^{\isom} &  {\hgm}(\O_{\cpt})\ar [r]^{\log} & K_2 \ar [r] & 0. }
\ee
Notably one has a surjective homomorphism of topological vector spaces
\be\label{eq:log-link}\widetilde{1+\fm_{K_1}} \to K_2.\ee
\ethm
\bp 
For any algebraically closed perfectoid field $K$ with $K^\flat=\cpt$ one has by \cite[Proposition 4.5.11]{fargues-fontaine}
$$\sgt(\O_K)\isom \sG(\O_{\cpt})$$
and a continuous surjection \cite[Proposition 4.5.14]{fargues-fontaine} 
$$\sG(\O_{\cpt}) \isom \sgt(\O_K) \onto K$$
which may be identified with the logarithm $\log_{\sG}$ of the formal group $\sG$.

This can be applied to the pair of fields $K_1,K_2$ and this provides the upper row of the diagram: applying the first isomorphism to $K_1$ one gets
$$\sgt(\O_{K_1})\isom \sG(\O_{\cpt})\isom \sgt(\O_{K_2}),$$
and composing it with the continuous surjection (provided by \cite[Proposition 4.5.14]{fargues-fontaine})
$$\sgt(\O_{K_2}) \onto K_2$$
one gets the top row of the diagram.

The lower row is obtained from the top row via the Artin-Hasse exponential mapping which provides an isomorphism of the formal groups $\sG\isom \hgm$ \cite[Example 4.4.7]{fargues-fontaine}. The  vertical arrows (isomorphisms) in the diagram are induced by the Artin-Hasse exponential. 
 
For $\hgm$, $\log_{\hgm}$ can be described explicitly in terms of the $p$-adic logarithm (\cite[Example 4.4.7]{fargues-fontaine})
\ep
 
 \brem 
In \iut\ a relationship between $p$-adic fields of the type \eqref{eq:log-link} is called a \emph{$\log$-link}. More precisely in \iut\ this would be a continuous homomorphism $\O_{K_1}^*/\text{Torsion}\mapright{\log} K_2$ (in the context of \iut, $K_1,K_2$ will be finite extensions of $\Q_p$). From the perspective of \cite{joshi-tempered}, Theorem~\ref{th:log-link}   is more natural and says that a version of Mochizuki's $\log$-link exists in the Arithmetic Teichmuller Theory of \cite{joshi-tempered} and this log-link does in fact relate two possibly distinct (perfectoid) valued field structures (especially if $K_1,K_2$ are not topologically  isomorphic algebraically closed perfectoid fields). In the parlance of (\cite[Page 6]{mochizuki-iut3}),   two, possibly distinct, ring structures (namely $K_1,K_2$) are linked by ``rendering additive'' the multiplicative structure. Here the multiplicative structure of $K_1$ is encapsulated in the multiplicative group $\widetilde{1+\fm_{K_1}}$ and the additive structure of $K_2$ is encapsulated in the topological group $(K_2,+)$ and a $\log$-link i.e. a morphism \eqref{eq:log-link} thus links the additive structure of $K_2$ encapsulated in the topological group $(K_2,+)$ to the multiplicative structure of $K_1$ embodied in  $\widetilde{1+\fm_{K_1}}$. In particular Theorem~\ref{th:log-link} shows log-links exist in the theory of \cite{joshi-tempered}.
 \erem

The following corollary of Theorem~\ref{th:log-link} is immediate.
\bcor\label{cor:log-link-in-arithmetic}
Let $K_1,K_2$ be two characteristic zero untilts of $\cpt$. Then 
\benumlab
\item one has an isomorphism
$$ \widetilde{1+\fm_{K_1}} \isom 1+\fm_{\cpt}\isom  \widetilde{1+\fm_{K_2}},$$
\item 
$K_1,K_2$ are log-linked via \eqref{eq:log-link} i.e. one has a surjection of topological vector spaces
$$\widetilde{1+\fm_{K_1}} \onto K_2,$$
where the arrow is the composite of $\widetilde{1+\fm_{K_1}}\mapright{\isom} 1+\fm_{\cpt} \mapright{\log} K_2$.
\item In general $K_1,K_2$ need not be topologically isomorphic.
\eenum
\ecor

\brem
Note that in \iut--especially in \cite{mochizuki-iut3}, Mochizuki works with $$\O_{\bQ_p}^{\times\mu}=\O_{\bQ_p}^*/\mu(\bQ_p)=(1+\fm_{ \bQ_p})/\mu_p(\bQ_p)$$
where $\mu(\bQ_p)$ (resp. $\mu_p(\bQ_p)$) is the subgroup of $\O_{\bQ_p}^*$ of roots of unity  (resp. subgroup of the $p$-power roots of unity) in $\bQ_p$, and especially with the exact sequence
$$1\to \mu_p(\bQ_p)\to 1+\fm_{ \bQ_p} \to \O_{\bQ_p}^{\times\mu}\to 1,$$
and considers ``different Kummer Theories'' (see \cite[Page 39--42]{mochizuki-gaussian}) i.e. different isomorphs of this sequence (and the associated Galois cohomology) obtained via Anabelian Reconstruction Algorithms of \topics, \iut,
while here and in \cite{joshi-tempered}, following \cite{fargues-fontaine}, I work with the exact sequence (dependent on $K$)
$$1\to T_p(\hgm)(\O_{K}) \tensor_{\Z_p}\Q_p \to \widetilde{1+\fm_{K}}=\widetilde{\hgm(\O_K)} \to K\to 0,$$
where $T_p(\hgm)(\O_K)$ is the ($p$-adic) Tate-module of $\hgm$ computed using $\O_K$ and $\hgm(\O_K)$ and allow $K$ to vary among algebraically closed perfectoid fields in general. 

Now let $$\bQ_{p;K}\subset K$$ be the algebraic closure of $\Q_p\subset K$ with the valuation induced from $K$. Then one has a Kummer sequence indexed by the algebraically closed perfectoid field $K$:
\be 1\to \mu_p(\bQ_{p;K})\to 1+\fm_{ \bQ_{p;K}} \to \O_{\bQ_{p;K}}^{\times\mu}\to 1.\ee
Especially if $K_1,K_2$ are not topologically isomorphic then one obtains two distinguishable Kummer sequences of the above type.
\erem

\iftoggle{arxiv}{
\bibliography{hoshi-bib,mochizuki-bib,uchida-bib,mochizuki-flowchart,../../master/master6.bib}
}
{
	\printbibliography
}

\Address
\end{document}